\documentclass[sigconf]{acmart}

\AtBeginDocument{%
  }

\renewcommand\footnotetextcopyrightpermission[1]{} 
\setcopyright{acmlicensed}
\copyrightyear{2018}
\acmYear{2018}
\acmDOI{XXXXXXX.XXXXXXX}

\acmConference[SFMES]{Simulation of Financial Markets and Economic Systems ICAIF 2024 Workshop}{14 November, 2024}{Brooklyn, NY}
\acmISBN{978-1-4503-XXXX-X/18/06}

\usepackage[linesnumbered,ruled,vlined]{algorithm2e} 

\begin{document}

\title{Deciding Bank Interest Rates}
\subtitle{A Major-Minor Impulse Control Mean-Field Game Perspective}

\author{Fan Chen}
\orcid{0000-0003-0082-7908}
\affiliation{%
    \institution{School of mathematical sciences, Shanghai Jiao Tong University}
    \city{Shanghai}
    \country{China}
}
\email{alexnwish@sjtu.edu.cn}

\author{Nicholas Martin}
\orcid{0000-0001-9266-1620}
\affiliation{%
  \city{London}
  \country{UK}}

\author{Po-Yu Chen}
\orcid{0000-0002-8568-8988}
\affiliation{%
  \institution{Imperial College London}
  \city{London}
  \country{UK}
}
\email{po-yu.chen11@imperial.ac.uk}

\author{Xiaozhen Wang}
\orcid{0009-0009-9286-0094}
\affiliation{%
 \institution{CEREMADE, Université Paris-Dauphine PSL}
 \city{Paris}
 \country{France}}
\email{xiaozhen.wang@dauphine.psl.eu}

\author{Zhenjie Ren}
\orcid{0000-0003-4656-4074}
\affiliation{%
  \institution{LaMME, Universit\'e Evry Paris-Saclay}
  \city{Paris}
  \country{France}}
\email{zhenjie.ren@univ-evry.fr}

\author{Francois Buet-Golfouse}
\orcid{0000-0002-2164-7087}
\affiliation{%
  \institution{AIML Global Markets, Barclays, University College London}
  \city{London}
  \country{UK}}
\email{francois.buetgolfouse@barclays.com}
\begin{abstract}
Deciding bank interest rates has been a long-standing challenge in finance. It is crucial to ensure that the selected rates balance market share and profitability. However, traditional approaches typically focus on the interest rate changes of individual banks, often neglecting the interactions with other banks in the market. This work proposes a novel framework that models the interest rate problem as a major-minor mean field game within the context of an interbank game. To incorporate the complex interactions between banks, we utilize mean-field theory and employ impulsive control to model the overhead in rate adjustments. Ultimately, we solve this optimal control problem using a new deep Q-network method, which iterates the parameterized action value functions for major and minor players and updates the networks in a Fictitious Play way. Our proposed algorithm converges, offering a solution that enables the analysis of strategies for major and minor players in the market under the Nash Equilibrium. 

\end{abstract}



\maketitle

\section{Introduction}
\label{sec:intro}

%
Effective competition among banks is an important safeguard for maintaining stability in the banking and financial sector 
 \cite{boyd2005theory}. Since the establishment of Bertrand competition model, bank can determine their share in the market by setting the number and price of financial products. 
Therefore, attracting more customers to bring in more revenue without sacrificing profitability by rationalising the interest rates on savings accounts has become the focus of attention in the inter-bank game \cite{van2013impact}, which leads to an optimal control problem 
\begin{quote}
\textit{"Find the optimal interest rate given the fact that there are other players in the market."}
\end{quote}
%
Most of the existing interest rate models focus only on the process of interest rate change of a single bank or institution. For example, the Vasicek model \cite{mamon2004three} is a short-rate model that describes the evolution of interest rates driven by only market risk. 
The Hull-White model \cite{hull1996using} is a more general formulation that parameterise the evolution of interest rate with respect to time.
However, although being frequently referred to by practitioners, none of the previous literature has addressed the problem from a multi-agent aspect, where multiple competitor banks aim to optimize the interest rate they offer and maximize profits over their market share.

In this work, we propose a novel framework that captures the interactions in the interbank game. This is the first framework that incorporates major-minor play theories to close the gap between well-known big banks and ordinary financial institutions (or the digital banks in the post-millennium era).
We then further exploit the mean field wheel to summarise the complex relationship of competition among banks. An impulsive control was also added to better incorporate the process of interest rate change into the proposed framework.
Based on the above framework, we specify the definition of the Nash equilibrium that guides the control strategies for the interbank game.

Ultimately, we solve the major-minor mean field game problem using the deep Q-network method. This novel approach numerically addresses major-minor mean field games, which are typically more complex due to the presence of common noise from the major player and the stochastic mean field measure flow. Moreover, our method applies to cases of continuous state space. To the best of our knowledge, there are no existing algorithms for this type of problem. The most relevant one is the algorithm proposed by Cui \cite{cui2024learning}, which is designed for finite state spaces. In contrast, our neural network-based approach can handle the general case and provide robust solutions.
\subsection{Contribution}
\label{subsec:contribution}
The contributions of this work are summarised below.
\begin{itemize}
    \item This is the first work to model the bank interest-rate problem as an interbank game. We establish a major-minor mean field game framework, mimicking the market role heterogeneity and incorporating impulsive controls to account for the adjustment costs in interest rate changes.
    \item We design a deep Q-network algorithm to solve the major-minor mean field game, allowing for continuous or large state and action spaces and overcoming the challenges of averaging neural networks inherent in Fictitious Play-based methods.
    \item The algorithm converges, providing robust solutions that present the strategies for major and minor players under Nash Equilibrium, supported by experimental results in Section \ref{sec:exp}.
\end{itemize}
\section{Related Work}
\subsection{Mean-Field Games}
The theory of Mean-Field Game (MFG) theory, introduced and developed by Lasry and Lions \cite{lasry2007mean} and Huang et al. \cite{ huang2006large}, incorporated the mean-field methodology widely used in fields of physics and chemistry into controlled games. In the regime of mean-field games, all atomic players are supposed to be homogeneous and symmetric and interact with each other. The number of players is assumed to be infinity, hence the interaction with other players can be characterized by the interaction with the distribution of players \cite{carmona2018probabilistic}. Recently, it has become a promising approach for modelling games with a large scale of interacting players 
\cite{tembine2020covid}. 
In finance and economics, extensive literature can be found for modelling market dynamics and systemic risks 
\cite{carmona2020applications}
, the construction of new price-impact models 
\cite{carmona2015probabilistic}
, the growth model in macro-economics 
\cite{achdou2022income}
, and the energy production and management 
\cite{elie2021mean}.
The pivotal objective of most research is the Nash Equilibrium, which is when each player in the game approaches its optimality against other players. For example, the existence analysis can be shown based on the contraction in fixed-point iterations \cite{carmona2018probabilistic}. 
We refer the reader to \cite{lauriere2022learning} 
for a more comprehensive review of applications and theoretical analysis for mean-field games in different settings.


\subsection{Major and Minor Players in MFG}

The major player was first introduced into a linear-quadratic infinite-horizon model by Huang \cite{huang2010large}. Distinguished with the assumption in classical MFG that each individual is symmetric and asymptotically negligible, there exist players who have a strong influence on the system. The major-minor MFG is surely an important case in practice. Pointed out as an example in \cite{chen2023risk}, in the price decision-making process, the decisions of innumerable minor companies are always affected by some major companies. Similar ideas are explored in analyzing interactions between large and small financial institutions by Bensoussan et al.\cite{bensoussan2016mean} and Carmona et al.\cite{carmona2014probabilistic}, and more applications can be found in \cite{dayanikli2023multi, carmona2022mean}. 


The major-minor MFG is considered a generalization of the common noise problem (see \cite{carmona2018probabilistic} for a basic framework). Since the common noise from the major player influences all minor players, the mean-field flow becomes stochastic, complicating both theoretical analysis and empirical implementation. The theory of major-minor MFG was further developed in the Linear-Quadratic-Gaussian (LQG) model 
\cite{firoozi2022lqg}. 
Other works proposed general frameworks under different settings, including formulating the Nash certainty equivalence principle 
\cite{huang2020linear}
, the stochastic Hamilton Jacobian Bellman 
\cite{bensoussan2016mean}
, and master equations 
\cite{cardaliaguet2020remarks}.
Another line of research focuses on the Stackelberg Equilibrium
\cite{guo2022optimization}
, in contrast with the Nash Equilibrium of all players. It is different from our settings, as in the Stackelberg Equilibrium, the major players (leaders) decide first and minor players (followers) respond to these decisions. 
%
%
\subsection{Algorithmic Approaches}
In this paper, we will focus on the computation of the Nash Equilibrium. Our approach is based on Reinforcement Learning (RL), a prevailing way with numerous applications (see survey \cite{wang2022deep}). 
For classical MFG, the basic methodology is alternating the estimation of population distribution and policy updates. 
%
A simple approach is to update the policy by the best response against other players \cite{guo2019learning}. 
There are several methods to stabilize the convergence process: 
the fictitious Play averages the past distributions or policies 
\cite{perrin2020fictitious}
; the online mirror descent and its variants accumulate the past Q-function 
\cite{perolat2021scaling}
; and other methods based on regularization ensure convergence, albeit with a bias 
\cite{cui2021approximately}.

When the state and action spaces are continuous, a tabular method will be infeasible and we should resort to the help of neural networks. However, the network parameterization of the policy or the Q-function faces the difficulty of taking averages, as emphasized in
\cite{lauriere2022scalable}. 
It is notable that in \cite{lauriere2022scalable}, the authors proposed the Munchausen Online Mirror Descent to overcome such difficulty, following the idea of Munchausen reparameterization \cite{vieillard2020munchausen}. In contrast, our proposed approach adopts a novel average technique by recognizing shallow neural networks as a function of measures as studied in \cite{chen2023entropic}, 
and our approach does not smooth the policy during iteration.

As for the Major-Minor player MFG, the most relevant paper to ours, Cui et al. designed a learning framework in discrete-time \cite{cui2024learning}. Fictitious play and projected mean-field algorithm were explored in \cite{cui2024learning}. Our Deep Q-Network (DQN) algorithm (follow the name in \cite{mnih2013playing}) is inspired by the thought of projection, nevertheless it is distinguished essentially from the one proposed in \cite{cui2024learning}: the average is upon the Q-function; it includes the learning procedure; the average of network outputs utilizes a measure perspective.
\subsection{Impulse Control}
Impulse controls are mathematically formulated to characterize the impulsive behaviours that widely exist in biology, medicine and economics. In impulse control problems, the agent can determine a sequence of controlled jumps to optimize its return. 
The impulse control problem has sparked considerable interest in applications, especially in resource economics 
\cite{alvarez2001stochastic}, portfolio optimization with transaction costs 
\cite{morton1995optimal}, inventory control and dividend control 
\cite{bai2010optimal}, and the modelling of exchange and interest rates 
\cite{mundaca1998optimal}.
Theoretically, the impulse control problems are naturally connected to the viscosity solution of Quasi-Variational Inequalities 
\cite{reisinger2020error}. 
Games of impulse control type are studied in 
\cite{azimzadeh2019zero,aid2020nonzero}
under various settings. We refer the readers to 
\cite{oksendal2019stochastic} 
for a comprehensive discussion and formulation of the impulse control problem.
\section{Framework}
\subsection{$N$-player Game}
Denote the time space as $\mathcal T$ and consider the continuous-time setting $\mathcal T = \mathbb R^+$. In the market, there are $M$ major banks, with each bank holding a proportion $(p_t^{(i)})_{i=1}^M$ of the entire market at time $t\in\mathcal T$. Additionally, there are $N$ minor banks holding proportions $(p_t^{(i)})_{i=M+1}^{M+N}$. These proportions satisfy the condition $\sum_{i=1}^{M+N}p_t^{(i)} = 1$. In the context of the mean-field game setting, the number of minor banks $N$ tends to infinity, and we assume that all minor banks are homogeneous while the major banks are not. Each bank will control its deposit interest rate $(r_t^i)_{i=1}^{M+N}$, which will influence its respective market proportion.
\paragraph{Proportion Change. } Let $(\kappa^{(i)})_{i=1}^{M+N}$ and $(\delta^{(i)})_{i=1}^{M+N}$ denote the escape rates and viscosity of clients, respectively. The dynamics of proportion change are characterized by:
\begin{equation}
\label{eq:density_transform}
\begin{aligned}
 \frac{dp_t^{(i)}}{dt} = \sum_{i'\neq i}\kappa^{(i')}&(r_t^{(i)}-r_t^{(i')}-\delta^{(i')})^+p_t^{(i)} \\
 &-{\sum_{i'\neq i} \kappa^{(i)}(r_t^{(i)}-r_t^{(i')}+\delta^{(i)})^-{p_t^{(i)}}}.
 \end{aligned}
\end{equation}
The viscosity parameters $(\delta^{(i)})_{i=1}^{M+N}$ can be understood as thresholds: once bank $i$ offers a higher rate than bank $i'$ to the extent that $r_t^{(i)}-r_t^{(i')}-\delta^{(i')}>0$, the clients of bank $i'$ changes their minds, causing the proportion to move from bank $i'$ to $i$.
It can also be observed that the larger the $\kappa^{i'}$ is, the more inclined the clients of bank $i'$ are to switch.
Furthermore, from \eqref{eq:density_transform}, it follows that the sum of all proportions remains constant over time, i.e., $\sum_{i=1}^{M+N}p_t^{(i)} = 1$.
\paragraph{Impulsive Control.} Banks' adjustments to deposit interest rates take the form of impulsive controls. Specifically, 
\begin{equation*}
    r_t^{(i)} = \sum_j \mathbf{1}_{\{t\ge\tau_j^{(i)}\}} \theta_j^{(i)}.
\end{equation*}
where $(\tau_j^{(i)},\theta_j^{(i)})_{j=1}^{\infty}$ is a sequence of stopping times and change magnitudes for bank $i$. In other words, each bank adjusts its interest rate at stopping times, together with deciding the most beneficial change magnitudes.
\paragraph{Accumulate P\&Ls.} Assume the total volume of bank deposits is fixed at $W$. At time $t$, the Profit \& Losses (P\&Ls) of bank $i$ from deposits are modeled by $Wp_t^{(i)}(l_t^{(i)}+r_t^c-r_t^{(i)})$, where $(l^i)_{i=1}^{M+N}$ represents the liquidity premium, $r_t^c$ denotes the Central Bank interest rate. Consequently, each bank faces a trade-off in setting its interest rate: it must avoid setting the rate too low, which would result in a loss of customers, or too high, which would result in excessive interest payments to depositors. The total discounted expected P\&L of the bank is given by:
{\footnotesize
\begin{equation*}
    J^{(i)} = \mathbb E\Bigg[\int_{\mathcal T} We^{-k^{(i)}t}p_t^{(i)}(l_t^{(i)}+r_t^c-r_t^{(i)}) dt -  \sum_j \mathbf{1}_{\{\tau_j^{(i)}\in\mathcal T\}} C^{(i)}(\theta_j^{(i)})\Bigg],
\end{equation*}}
where $C^{(i)}(\theta_j^{(i)})>0$ represents the cost incurred by banks in adjusting interest rates, and $e^{k^{(i)}t}$ is the discounting factor. As a result, banks cannot make adjustments continuously or without careful consideration.
\paragraph{Central Bank Rate.}   
In our model, the Central Bank rate is described by the jump process
\begin{equation}
\label{eq:CB_rate_dynamic}
    dr_t^c = \int_{\mathbb{R}} \alpha(t, r_{t-}^c, z) N(dt, dz),
\end{equation}
where $ N(dt, dz) $ is the Poisson random measure (refer to \cite{oksendal2019stochastic} for its definition). The integral variable $z$ represents the size of a jump that happens in the Poisson process, $ \alpha(t, r_{t-}^c, z) $ characterizes the magnitude of these jumps, and $ r_{t-} $ denotes the left limit of the rate process. The randomness of the central bank rate serves as an external source of noise, affecting all players simultaneously. In numerical experiments, we discrete the time space and treat the central bank rate dynamics as random transitions between states.
\subsection{Mean-Field Regime}
All minor banks are assumed to be homogeneous, with parameters $\delta^{(i)} = \delta$, $\kappa^{(i)} =\kappa$ for $i=M+1,...,M+N$.
We assume there are infinitely many minor banks. The proportion change dynamics in \eqref{eq:density_transform} are modified to
\begin{equation}
\label{eq:density_transform_modified}
{
    \begin{aligned}
        \frac{d p_t^{(i)}}{dt}
        = &\sum_{i'=1}^M\kappa^{(i')}(r_t^{(i)}-r_t^{(i')}-\delta^{(i')})^+ p_t^{(i')} \\
        &\quad + \frac{1}{N}\sum_{i'=M+1}^{M+N}\kappa(r_t^{(i)}-r_t-\delta)^+ p_t \\
        &-{\sum_{i'=1}^M \kappa^{(i)}(r_t^{(i)}-r_t^{(i')}+\delta^{(i)})^-{ p_t^{(i)}}} \\
        &\quad - \frac{1}{N}\sum_{i'=M+1}^{M+N} \kappa(r_t^{(i)}-r_t+\delta^{(i)})^- p_t^{(i)},
    \end{aligned}
}%
\end{equation}
where we set a multiplier $1/N$ to rescale the proportion exchanges between each pair of major and minor banks. 
Rescaling the proportion of minor banks as $\bar p_t^{(i)} = Np_t^{(i)}$ for $i=M+1,...,M+N$ and abusing the notation $\bar p_t^{(i)} = p_t^{(i)}$ for $i=1,...,M$, we are ready to rewrite \eqref{eq:density_transform_modified} as
\begin{equation}
\label{eq:density_transform_modified_integral_major}
{
    \begin{aligned}
     \frac{d\bar p_t^{(i)}}{dt} = &\sum_{i'=1}^M\kappa^{(i')}(r_t^{(i)}-r_t^{(i')}-\delta^{(i')})^+\bar p_t^{(i')} \\
     &\quad + \int \kappa(r_t^{(i)}-r_t-\delta)^+\bar p_t\mu_t^N(d\bar p_t,dr_t)  \\
     &-{\sum_{i'=1}^M \kappa^{(i)}(r_t^{(i)}-r_t^{(i')}+\delta^{(i)})^-{\bar p_t^{(i)}}}\\
     &\quad -\int \kappa(r_t^{(i)}-r_t+\delta^{(i)})^-\bar p_t^{(i)}\mu_t^N(d\bar p_t,dr_t),
    \end{aligned}
}%
\end{equation}
for all banks $i=1,2,...,M+N$, where $\mu_t^N = \sum_{i=M+1}^{M+N} \delta_{(\bar p_t^{(i)},r_t^{(i)})}$. The rescaling of dynamics and proportions ensures a proper definition of the distribution $\mu_t^N$ as $N$ becomes large. Under the mean-field framework, the dynamics of representative \textit{minor players} are given by
\begin{equation}
\label{eq:density_transform_modified_integral_minor_rep}
{\
    \begin{aligned}
    \frac{d\bar p_t}{dt} 
    = &\sum_{i'=1}^M\kappa^{(i')}(r_t^{(i)}-r_t^{(i')}-\delta^{(i')})^+\bar p_t^{(i')} \\
    &\quad + \int \kappa(r_t^{(i)}-r_t-\delta)^+ p_t\mu_t(d p_t,dr_t)  \\
    &-{\sum_{i'=1}^M \kappa^{(i)}(r_t^{(i)}-r_t^{(i')}+\delta)^-{\bar p_t}}\\
    &\quad -\int \kappa(r_t^{(i)}-r_t+\delta)^-\bar p_t\mu_t(d p_t,dr_t), \\
    \eqcolon &~~b(\mathbf{\bar p_t^M},\mathbf{r_t^M},\bar p_t,r_t,\mu_t),
    \end{aligned}
}%
\end{equation}

where $\mathbf{\bar p_t^M}$ is the collection of all major players' proportions, and similarly, $\mathbf{r_t^M}$ is the collection of their interest rates.  Here, $\bar p_t$ is regarded as a stochastic process with mean field flow $\mu_t = \text{Law}(\bar p_t,r_t)$ at each time $t$. On the other hand, the dynamics for the \textit{major players} can be expressed as:
\begin{equation}
\label{eq:density_transform_modified_integral_major_MF}
{
    \begin{aligned}
     \frac{d\bar p_t^{(i)}}{dt} = &\sum_{i'=1}^M\kappa^{(i')}(r_t^{(i)}-r_t^{(i')}-\delta^{(i')})^+\bar p_t^{(i')} \\
     &\quad + \int \kappa(r_t^{(i)}-r_t-\delta)^+\bar p_t\mu_t(d\bar p_t,dr_t)  \\
     &-{\sum_{i'=1}^M \kappa^{(i)}(r_t^{(i)}-r_t^{(i')}+\delta^{(i)})^-{\bar p_t^{(i)}}}\\
     &\quad -\int \kappa(r_t^{(i)}-r_t+\delta^{(i)})^-\bar p_t^{(i)}\mu_t(d\bar p_t,dr_t), \\
     \eqcolon &~~b^{(i)}(\mathbf{\bar p_t^M},\mathbf{r_t^M},\mu_t)
    \end{aligned}
}%
\end{equation}
and again it satisfies that $\sum_{i=1}^M \bar p_t^{(i)} + \int \bar p_t \mu_t(d\bar p_t,dr_t) = 1$ for all $t\in\mathcal T$. Consequently, the objective functions of the banks take the form
\begin{equation*}
{
    \begin{aligned}
        J^{(i)}&(\mathbf{\tau}^{(i)},\mathbf \theta^{(i)}) :=\\
        &\mathbb E\Bigg[\int_{\mathcal T} We^{-k^{(i)}t}p_t^{(i)}(l_t^{(i)}+r_t^c-r_t^{(i)}) dt -  \sum_j \mathbf{1}_{\{\tau_j^{(i)}\in\mathcal T\}} C^{(i)}(\theta_j^{(i)})\Bigg], \\
        &J(\mathbf \tau,\mathbf \theta) :=\mathbb E\Bigg[\int_{\mathcal T} We^{-kt}\bar p_t(l_t+r_t^c-r_t) dt -  \sum_j \mathbf{1}_{\{\tau_j\in\mathcal T\}} C(\theta_j)\Bigg], \\
    \end{aligned}
}
\end{equation*}
where $\tau^{(i)} = \{\tau_j^{(i)}\}_{j=1}^\infty$ (resp. $\tau = \{\tau_j\}_{j=1}^\infty$) and $\theta^{(i)} = \{\theta_j^{(i)}\}_{j=1}^\infty$ (resp. $\theta = \{\theta_j\}_{j=1}^\infty$) are the controls decided by the major players (resp. 
representative minor player).
\subsection{Discrete-time Model}
For simplicity, we now consider the case where there is one major bank and multiple minor banks. The time space is $\mathcal{T}=[T]$. Initially, the state of the major player is fixed as $x_0^0$ and the distribution of the minor players is set to $\mu_0$, which is assumed to be supported on a compact set. At time $t\in [T]$, the states of the major bank are represented as $x_t^0 = (p_t^0,r_t^0)$, while for the representative minor bank, they are $x_t = (p_t,r_t)$. Along with the states, the central bank rate is denoted as $r_t^c$, and the mean-field flow $\mu_t$. In the discrete-time setting, impulse controls are treated similarly to regular controls. For any $t\in [T]$, the major and minor players will control $(u_t^0,u_t)$ to be their new interest rates, i.e., $r_{t+1}^0 = u_t^0$ and $r_{t+1} = u_t$. Subsequently, the proportions transfer according to these new interest rates. 

Let $\mathcal X$ and $\mathcal U$ be the state and action space for minor players,  $\mathcal X^0$ and $\mathcal U^0$ for the major player, and let $\mathcal R$ be the space of central bank rate. Let $K^0(x^0,u^0,\mu):\mathcal X^0\times \mathcal U^0\times \mathcal P(\mathcal X)\to \mathcal X^0$ and $K(x^0,u^0,x,u,\mu):\mathcal X^0\times \mathcal U^0\times \mathcal X\times \mathcal U\times \mathcal P(\mathcal X)\to \mathcal X$ be deterministic transition functions such that 
\begin{equation*}
{
\begin{aligned}
    &K_p^0(x^0,u^0,\mu) =  p^0 + b^0(x^0,u^0,\mu)\Delta t, \ \
    K_r^0(x^0,u^0,\mu) =   u^0,\\
    &K_p(x^0,u^0,x,u,\mu) = x +  b(x^0,u^0,x,u,\mu)\Delta t, \ \  
    K_r(x^0,u^0,x,u,\mu) =  u,\\
    \end{aligned}
}
\end{equation*}
with $K^0 = (K^0_p,K^0_r)$ and $K = (K_p,K_r)$. The drifts $b^0$ and $b$ are defined in \eqref{eq:density_transform_modified_integral_minor_rep} and \eqref{eq:density_transform_modified_integral_major_MF} for the case $M=1$. Let $P^c(r^c):\mathcal R\to P(\mathcal R)$ be the transition probability for the central bank rate. The states of players evolve as $x_{t+1}^0  = K^0\left(x_t^0, u_t^0, \mu_t\right)$ for the major player and $x_{t+1}  = P\left(x_t^0, u_t^0, x_t,u_t, \mu_t\right)$ for the minor players. The central bank rate transits as $r^c_{t+1} \sim P^c(r_t^c)$. Let the control function of major player be $u_t^0(x^0,r^c,\mu):\mathcal X^0\times \mathcal R\times \mathcal P(\mathcal X)\to \mathcal U^0$, and for minor players $u_t(x^0,x,r^c,\mu):\mathcal X^0\times \mathcal X\times \mathcal R\times \mathcal P(\mathcal X)\to \mathcal U$. Given a sequential control functions $\mathbf u =( u_t)_{t\in\mathcal T}$, we define the operator 
\begin{equation*}
{
    T_t^{\mathbf u}\left(x_t^0, u_t^0, r_t^c,\mu_t\right):=\int_{\mathcal X } K\Big(x_t^0, u_t^0,x,  u_t(x_t^0,x, r_t^c, \mu_t), \mu\Big) \mu_t(d x)
}%
\end{equation*}
to represent the mean field transition $\mu_t$ to $\mu_{t+1}$
, if all the minor plays apply the control function ${\mathbf u}_t$ and the major player takes state $x_t^0$ and action $u_t^0$.

If we also define a sequence of major control functions as $\mathbf u^0 = (\mathbf u_t^0)_{t\in\mathcal T}$, the objective functions are defined as follows:
\begin{equation*}
{
    \begin{aligned}
    &J^0\left(\mathbf u^0, \mathbf u\right)=\mathbb{E}\left[\sum_{t \in [T]} \gamma^tR_t^0\left(x_t^0, u_t^0(x_t^0,r_t^c, \mu_t), r_t^c,\mu_t\right)\right], \\
    &J\left(\mathbf u^0, \mathbf u\right)=\mathbb{E}\left[\sum_{t \in [T]}\gamma^t R_t\left(x_t, u_t(x_t^0,x_t, r_t^c, \mu_t), r_t^c, \mu_t\right)\right], \\
    \end{aligned}
}%
\end{equation*}
where $\gamma=\exp(-k^0\Delta t)=\exp(-k\Delta t)$ , and the running rewards are given by
\begin{equation*}
{
    \begin{aligned}
        &R_t^0\left(x_t^0, u_t^0, r_t^c,\mu_t\right) = W  p_{t}^0(l^0+r_t^c-u_t)-C^0(u_t^0-r_t^0),\\
        &R_t\left(x_t, u_t, r_t^c, \mu_t\right) = Wp_{t}(l+r_t^c-u_t)-C(u_t-r_t).
    \end{aligned}
}%
\end{equation*}
\subsection{Nash Equilibrium}
To define the Nash Equilibrium, we introduce an individual minor player who follows the control functions $\mathbf {\hat u}$, while other minor players continue to follow the controls $\mathbf u$. Note that this single minor player does not influence the evolution of other players. The objective function for this player is defined as:
\begin{equation*}
\begin{aligned}
J\left(\mathbf u^0, \mathbf u,\mathbf {\hat u}\right)=\mathbb{E}\left[\sum_{t \in [T]} \gamma^t R_t\left(x_t^0, u_t^0, \hat x_t,\hat u_t ,\mu_t\right)\right]. \\
\end{aligned}
\end{equation*}

\begin{definition}
\label{def:nash}
    In our model, a control profile $(\mathbf u_*^0,\mathbf u_*)$ is called a Nash Equilibrium if for any alternative policies $\mathbf u^0$ and $\mathbf {\hat u}$, it satisfies the conditions
\begin{equation*}
    \begin{aligned}
        &J^0\left(\mathbf u_*^0,\mathbf u_*\right) \geq J^0\left(\mathbf u^0,\mathbf u_*\right),\\
        &J(\mathbf u_*^0,\mathbf u_*,\mathbf u_*) \geq J\left(\mathbf u_*^0,\mathbf u_*,\mathbf {\hat u}\right).
    \end{aligned}
\end{equation*}
\end{definition}
\noindent
The Nash Equilibrium conditions indicate that the controls are optimal against those of other players. For the minor players, the controls are optimal even against other minor players. Therefore, according to Definition \ref{def:nash}, we aim to solve a fixed point problem. 


\section{Algorithm}
Now, we introduce our Deep Q-Network (DQN) algorithm, based on the parametrization of the action-value function (Q function). The algorithm aims to ensure that Bellman's equation holds or approximately holds. For example for the major player, the action value function, given the conjectured policies of the other players $\mathbf u$, should satisfy
\begin{equation}
\label{eq:Bellman_eq}
{
    \begin{aligned}
    Q_{\mathbf u, \mathbf u^0}^0&(t, x^0, u^0, r^c, \mu)=\gamma^{t}  R_t^0(x^0, u^0, r^c, \mu) + \max _{\tilde u^{0}\in\mathcal U^0} \int_{\mathcal R} P^c(d\tilde r^c| r^c)\\
    \cdot& Q_{\mathbf u, \mathbf u^0}^0\big(t+1, K^0(x^0,u^0,\mu), \tilde u^{0}, \tilde r^c, T_t^{\mathbf u}(x^0, u^0, r^c, \mu)\big) .   
    \end{aligned}
}
\end{equation}
Note that this equation is equivalent to the one where \(\gamma\) appears before the integral, as commonly used in the reinforcement learning literature (see \cite{wiering2012reinforcement}). We can consider a discounted Q function to obtain this formulation.

\paragraph{Projection of Mean-Field Measure.} To parametrize the Q function by a neural network, we approximate the measure $\mu\in\mathcal P(\mathcal X)$ using finite-dimensional vectors. It is necessary since the measure object lies in an infinite dimensional space and we should consider the entire measure as an input. This process is called the projection of a mean-field measure. It is feasible given the assumption that $\mu_0$ is compactly supported, we can estimate a compact subset $\mathcal X^c$ such that $\text{supp}(\mu_t)\subset \mathcal X^c$ for all $t\in \mathcal T$, since both $b^0$ and $b$ are bounded. Let $\mathcal X^F = (x^{(i)})_{i=1}^{N_\epsilon}$ be a $\epsilon-$net for the compact set $\mathcal X^c$. Elements of $\mu\in \mathcal P(\mathcal X^F)$ can be expressed as vectors $N_\epsilon$-length of length $N_\epsilon$. We define the projection operator $\mathcal A:\mathcal P(\mathcal X^c)\to \mathcal P(\mathcal X^F)$ such that $\mathcal A(\mu_t)\in \mathcal P(\mathcal X^F)$ is a good approximation of $\mu_t\in \mathcal P(\mathcal X^c)$.
Then, we modify Bellman's equation \eqref{eq:Bellman_eq} as follows
\begin{equation}
\label{eq:Bellman_eq_modified}
{
    \begin{aligned}
    &\hat Q_{\mathbf u, \mathbf u^0}^0(t, x^0, u^0, r^c, \mu)=\hat Q_{\mathbf u, \mathbf u^0}^0(t, x^0, u^0, r^c, \mathcal A(\mu))\\
    &=\gamma^t R_t^0(x^0, u^0, r^c, \mathcal A(\mu))+\max _{\tilde u^{0}\in \mathcal U^0}\int_{\mathcal R} P^c(d\tilde r^c| r^c) \\
    &\cdot  \hat Q_{\mathbf u, \mathbf u^0}^0\Big(t+1, K^0(x^0,u^0,\mathcal A(\mu)), \tilde u^{0}, \tilde r^c,\mathcal A\big(T_t^{\mathbf u}(x^0, u^0,  \tilde r^c, \mathcal A(\mu))\big)\Big).
    \end{aligned}
}%
\end{equation}
Based on the approximated Bellman's equation \eqref{eq:Bellman_eq_modified}, we represent $\hat Q_{\mathbf u, \mathbf u^0}^0(t, x^0, u^0, r^c, \mathcal A(\mu))$ as a neural network in the form $\hat Q(\theta;\omega)$, where the inputs are $\theta^0 = \big(t,x^0,u^0,r^c,(\mu_t^{(i)})_{i=1}^{N_\epsilon}\big)$, with $\sum_{i=1}^{N_\epsilon} \mu_t^{(i)}= 1$ and $\mu_t^{(i)}$ representing the mass at point $x_t^i$. We use optimization algorithms for neural networks to update the load parameters $\omega$. Similarly, we derive the approximated equality for the Q-function of the representative minor player as
\begin{equation}
\label{eq:Bellman_eq_modified_minor}
{
    \begin{aligned}
    &\hat Q_{\mathbf u, \mathbf u^0}(t, x^0, x,u,r^c, \mu)
    =\hat Q_{\mathbf u, \mathbf u^0}(t, x^0,x,u, r^c, \mathcal A(\mu))\\
    &=\gamma^t R_t(x, u, r^c, \mathcal A(\mu))+\max _{\tilde u\in \mathcal U}\int_{\mathcal R} P^c(d\tilde r^c| r^c) \\
    &\cdot  \hat Q_{\mathbf u, \mathbf u^0}\Big(t+1, 
    K^0(\cdot), 
    K(\cdot),
    \tilde u, \tilde r^c,\mathcal A\big(T_t^{\mathbf u}(x^0, \mathbf u^0, r^c, \mathcal A(\mu))\big)\Big)
    \end{aligned}
}%
\end{equation}
where {\small $K^0(\cdot) \coloneq K^0(x^0,\mathbf u^0,\mathcal A(\mu))$} and {\small $K(\cdot) \coloneq K(x^0,u^0,x, u,\mathcal A(\mu))$}.
We use the function $\mathcal S \hat Q^0(\cdot)$ and $\mathcal S\hat  Q(\cdot)$ to denote the right-hand side of the equations \eqref{eq:Bellman_eq_modified} and \eqref{eq:Bellman_eq_modified_minor}, respectively. 

\begin{figure}
    \centering
    \includegraphics[scale=0.15]{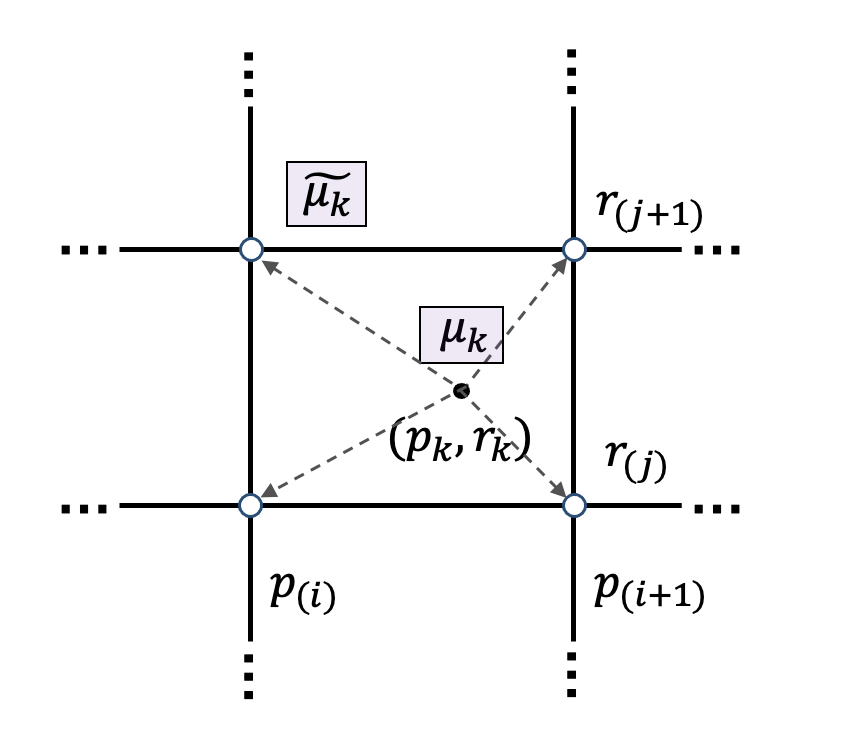}
    \caption{Diagram of constructing projection $\mathcal A$: $\mu_k\to \tilde \mu_k$}
    \label{fig:sample_projection}
\end{figure}

\paragraph{Choice of Projection Operator $\mathcal A$}
In fact, in training, we are only concerned with projecting probability measure supported on finite points to $\mathcal P (\mathcal X^F)$. Denote the original probability measure $\mu$ with supports $(p_k,r_k)_{k\in [K]}$ where $\mu_k = \mu((p_k,r_k))$ and $\sum_{k\in[K]} \mu_k = 1$. For each $k\in [K]$, the corresponding support $(p_k,r_k)$ falls into a sub-rectangle, i.e., $(p_k,r_k)\in (p_{(i)},p_{(i+1)})\times (r_{(j)},r_{(j+1)})$. We define a new measure $\tilde \mu^k$ supported on the vertices of this sub-rectangle (see Figure \ref{fig:sample_projection}) as
\begin{equation*}
    \tilde \mu_k(p,r) = \left\{\begin{array}{ll}
         (p_{(i+1)}-p_k)(r_{(j+1)}-r_k)\mu_k, &\text{if}\ (p,r)=(p_{(i)},r_{(j)}); \\
          (p_{(i+1)}-p_k)(r_k - r_{(j)})\mu_k, &\text{if}\ (p,r)=(p_{(i)},r_{(j+1)}); \\
          (p_k-p_{(i)})(r_{(j+1)}-r_k)\mu_k, &\text{if}\ (p,r)=(p_{(i+1)},r_{(j)}); \\
          (p_k-p_{(i)})(r_k - r_{(j)})\mu_k, &\text{if}\ (p,r)=(p_{(i+1)},r_{(j+1)}); \\
          0, &\text{otherwise}.
    \end{array}\right.
\end{equation*}

As an approximation of the original finite supported distribution $\mu$, we define $\mathcal A(\mu) = \sum_{k\in [K]} \tilde \mu_k\in\mathcal P(\mathcal X^F)$.

\paragraph{Fictitious Play and Average Neural Networks}
To stabilize the convergence process, we use a moving average of the past Q functions. Denoting $Q_n^0$ and $Q_n$ as the approximated Q-functions at the $n$-th iteration, we derive the optimal control functions $u_n^0$ and $u_n$ from these Q functions. We then evaluate the Q-functions according to the Bellman Equations \eqref{eq:Bellman_eq_modified} and \eqref{eq:Bellman_eq_modified_minor} and the strategies of the opposite players. Let $\hat Q_{n}^0$, $\hat Q_{n}$ represent the evaluated Q functions. We update the Q-function as $Q_{n+1} \leftarrow \frac{n}{n+1} Q_{n} + \frac{1}{n+1} \hat Q_{n} $, $Q_{n+1}^0 \leftarrow \frac{n}{n+1} Q_{n}^0 + \frac{1}{n+1} \hat Q_{n}^0 $ and pass to the next step (see Algorithm \ref{alg:DQN}).

The challenge lies in averaging neural networks due to their non-linear nature, as mentioned in \cite{lauriere2022scalable}. In our experiments, we use fully connected neural networks with one hidden layer as our architecture. The parameters of neural networks can be interpreted as measures, as described in \cite{chen2023entropic, hu2021mean}. It can be represented as
\begin{equation}
\label{eq:measure_landscape_intro}
\frac{1}{L} \sum_{l=1}^L \beta_{l} \varphi\left(\alpha_{ l} \cdot \theta+\gamma_{l}\right)=\int_{\mathbb{R}^{d+2}} \beta \varphi(\alpha \cdot \theta+\gamma) m^L(d \beta d \alpha d \gamma)   
\end{equation}
where we denote the input $z\in\mathbb R^d$, parameters as multiplier $\alpha$, bias $\gamma$ and weight $\beta$, and $m^L=\frac{1}{L} \sum_{l=1}^L \delta_{\left(\beta_{l}, \alpha_{l}, \gamma_{l}\right)}.$
For simplicity, denote $x=(\beta, \alpha, \gamma) \in \mathbb{R}^{d+2}$ and $\hat{\varphi}(x, z)=\beta \varphi(\alpha \cdot z+\gamma)$, so the neural network output can be represented as $\hat Q(\theta;m)= \mathbb{E}^{X\sim m}[\hat{\varphi}(X, z)]$. Thus, when the width $L$ becomes large, we can understand the load parameters as a measure. This allows us to transfer the problem from averaging the outputs of neural networks to averaging the measures, which is more feasible in implementation. To be precise,
we parameterized $\hat Q(\theta;m_{n}) = \mathbb{E}^{X\sim m_n}[\hat{\varphi}(X, \theta)]$, $\hat Q(\theta;\tilde m_n) = \mathbb{E}^{X\sim \tilde m_{n}}[\hat{\varphi}(X, \theta)]$, the update of the moving average can be expressed as 
\begin{equation}
\label{eq:ficititous_play_process}
\begin{aligned}
    &\hat Q(z;m_{n+1}) = \frac{n}{n+1}\hat Q(z;m_{n}) + \frac{1}{n+1} \hat Q(z;\tilde m_n) \\
    &\quad \quad \Leftrightarrow  m_{n+1} = \frac{n}{n+1} m_{n} + \frac{1}{n+1} \tilde m_{n}.    
\end{aligned}
\end{equation}
The weighted average of the measures can be achieved by combining neurons through random selection. 
\begin{algorithm}
\caption{Deep Q-Network (DQN) Algorithm}
\label{alg:DQN}
\KwIn{Approximate $\epsilon$-net $\mathcal X^F$, transition functions $K,K^0$ and $P^c$, drift function $b^0,b$, the reward function $R^0,R$, iteration number of outer iteration $N$, iteration number of inner iteration for network update $M$, batch size $B$;}
Initialize network $\hat Q^0(\theta^0;\omega_{0,0}^0)$ and $\hat Q(\theta;\omega_{0,0})$;\\
\For{$n=0,1,..N-1$}{
Set $u_t^0(x^0,r^c,\mu) = \text{argmax}_{u^0\in \mathcal U^0} \hat Q^0(t,x^0,u^0,r^c,\mu;\omega_{n,0}^0)$, $u_t(x^0,x,r^c,\mu) = \text{argmax}_{u\in \mathcal U} \hat Q^0(t,x^0,x,u,r^c,\mu;\omega_{n,0})$ for $t\in\mathcal T$, $x\in\mathcal X$, $x^0\in \mathcal X^0$, $r^c\in \mathcal R$, $\mu\in \mathcal P(\mathcal X^F)$ ;\\
\For{$m=0,1,..M-1$}{
Sample a batch of inputs as $\big\{(t_i,x^0_i,u_i^0,x_i,u_i,r^c_i,(\mu^{(j)}_i)_{j=1}^{N_\epsilon}\big\}_{i=1}^{B}$;\\
Update $\omega_{n,m}^0\mapsto  \omega_{n,m+1}^0$ by optimizers based on loss:
$$
L^0(\omega^0) = \frac{1}{B}\sum_{i=1}^{B}\big|\hat Q^0(\theta_i^0;\omega^0) - \mathcal S  \hat Q^0(\theta_i^0;\omega^0) \big|^2;
$$\\
Update $\omega_{n,m}\mapsto  \omega_{n, m+1}$ by optimizers based on loss:
$$
L(\omega) = \frac{1}{B}\sum_{i=1}^{B}\big|\hat Q(\theta_i;\omega) - \mathcal S  \hat Q(\theta_i;\omega) \big|^2;
$$}
Update $Q^0,Q$ in a fictitious play way, as in \eqref{eq:ficititous_play_process}:
$$
\begin{aligned}
&\hat Q^0(\theta^0;\omega_{n+1,0}^0) = \frac{n}{n+1} \hat Q^0(\theta^0;\omega_{n,0}^0) + \frac{1}{n+1}\hat Q^0(\theta^0; \omega_{n,M}^0); \\
&\hat Q(\theta;\omega_{n+1,0}) = \frac{n}{n+1}\hat Q(\theta;\omega_{n,0}) + \frac{1}{n+1}\hat Q(\theta; \omega_{n,M}) ;  
\end{aligned}
$$
}
Define the control functions:
\begin{equation*}
    \begin{aligned}
        &u_t^0(x^0,r^c,\mu) = u_t^0(x^0,r^c,\mathcal A(\mu)) = \\
        &\quad \quad \quad \quad\quad \quad\text{argmax}_{u^0\in \mathcal U^0} \hat Q^0(t,x^0,u^0,r^c,\mu;\omega_{N,0}^0),\\
&u_t(x^0,x,R,\mu) = u_t(x^0,x,r^c,\mathcal A(\mu))   = \\
&\quad \quad \quad \quad\quad \quad\text{argmax}_{u\in \mathcal U} \hat Q^0(t,x^0,x,u,r^c,\mu;\omega_{N,0}),
\end{aligned}
\end{equation*}
for any $t\in\mathcal T$, $x\in\mathcal X$, $x^0\in \mathcal X^0$, $r^c\in \mathcal R$, $\mu\in \mathcal P(\mathcal X)$;
\\ 
\KwOut{The control functions $u_t^0$ and $u_t$.}
\end{algorithm}

\begin{figure*}[htbp!]
\begin{minipage}{0.48\linewidth}
\centering
\includegraphics[width=\linewidth,height=4cm]{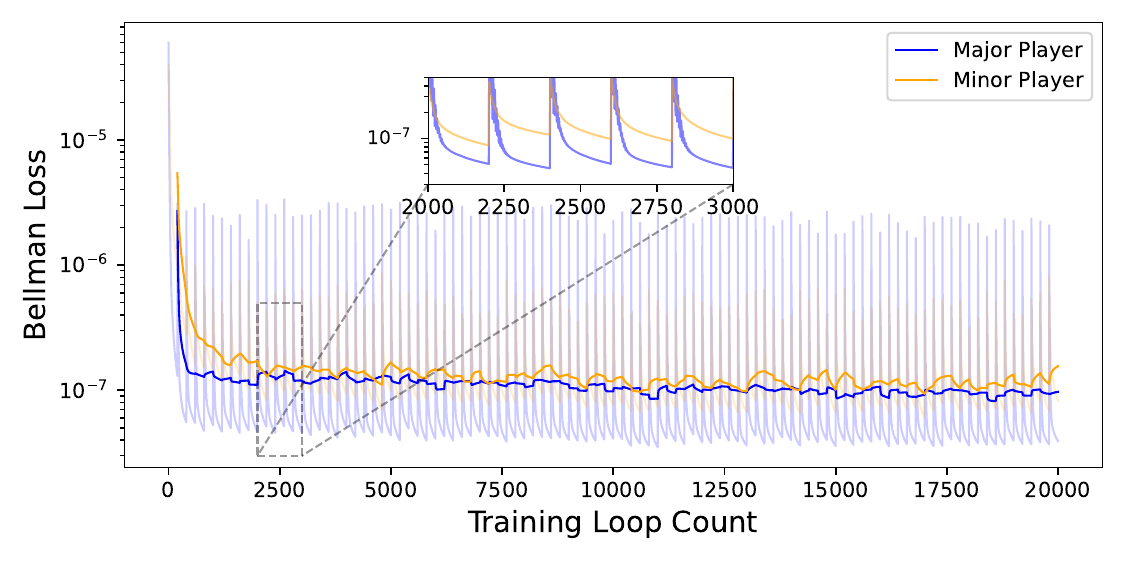}
\captionof{figure}{Average losses of 10 individual
trainings (shadowed) and its \(200\) rolling average (bold).}
\label{fig:error-gap}
\end{minipage}
\hfill
\begin{minipage}{0.48\linewidth}
\centering
\includegraphics[width=\linewidth,height=4cm]{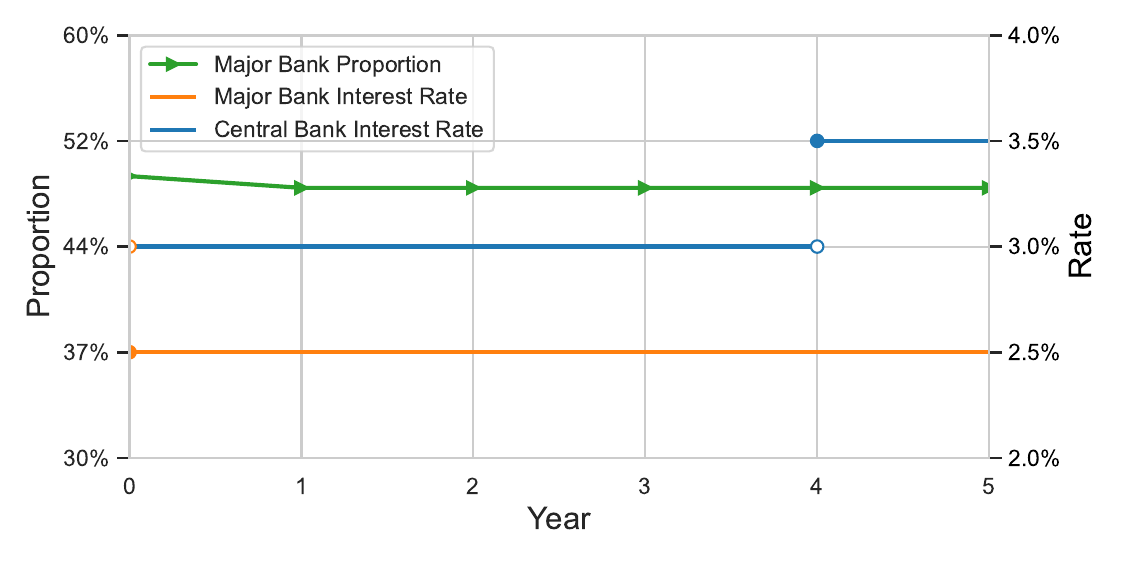}
\captionof{figure}{The evolution of major player's state and central bank rate. }
\label{fig:dynamic_evolution}
\end{minipage}
\end{figure*}

\begin{figure*}[htbp!]
    \centering
    \includegraphics[width=\textwidth]{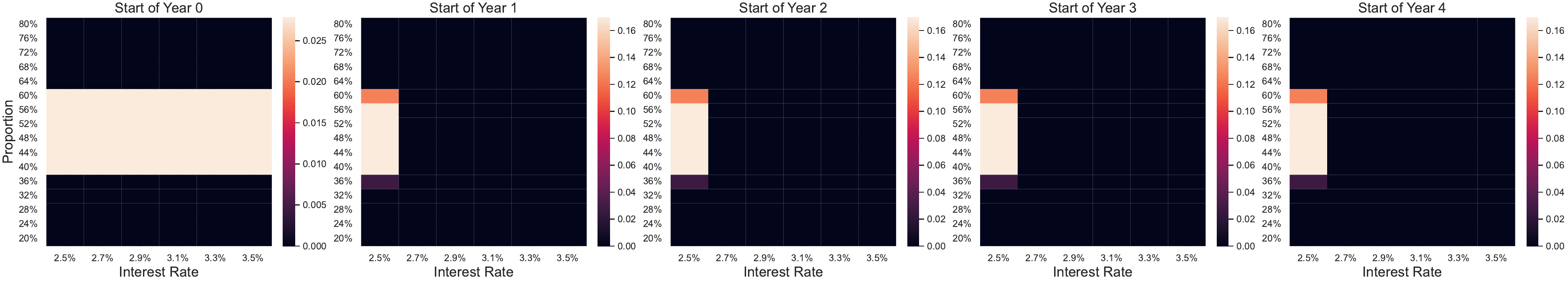}
    \caption{The evolution of mean field measure (minor players).}
    \label{fig:mean_field_evolution}
\end{figure*}
\section{Experiments}
\label{sec:exp}
\subsection{Setup}
In our experiment, we assume that the time horizon of the banks is $T=5$, with a decision period $\Delta t =1$. For clients' preference, the escape rate and the viscosity parameter of banks, we set $(\kappa^0, \kappa)=(5,5)$ and $(\delta^0,\delta)=(0.1\%,0.1\%)$. The range of acceptable interest rate selection $\mathcal U^0 = \mathcal U = [2.5\%,3.5\%]$, the state space $\mathcal X^0 = \mathcal X = [20\%,80\%] \times [2.5\%,3.5\%]$. Given these parameters and the definition of $b^0,b$ (formulated in \eqref{eq:density_transform_modified_integral_minor_rep} and \eqref{eq:density_transform_modified_integral_major_MF}), we have
\begin{equation}
\label{eq:b_bound}
    |b^0|,|b|\le \overline \kappa (r_{max}-r_{min} - \underline{\delta}) = 4.5\%.
\end{equation}
This implies that the proportion of major and minor players will change by at most $2\%$ due to the interest rate competition. Initially, for the state of the major player, we set $\bar p^0=50\%$, $r^0=3\%$, while for the minor players, we set their distribution to be uniform $\mu_0=U([40\%.60\%]\times [2.5\%,3.5\%])$. Given the bound $b^0$, $b$ \eqref{eq:b_bound}, we deduce that the compact support of $\mu_t$ will not escape from $[20\%.80\%]\times [2.5\%,3.5\%]$, that is, $\mathcal X^c = \mathcal X$. Then, we evenly divide the space $\mathcal X^c$ into $12\times 5$ sub-rectangles. Let $\mathcal X^F = (p_{(i)},r_{(j)})_{i\in[15],j\in[5]}$ with $p_{(i)} = 20\% + 4\%\cdot i$ and $r_{(j)}=  2.5\% + 0.2\%\cdot j$. Throughout our paper, we set the discounting $\gamma = 0.9$ and the total deposit volume $W=1$ without loss of generality. The liquidity premiums are set to $l^0=0$, $l=0.1\%$. Regarding the cost incurred from interest rate changes, we set the cost function $C^0(\Delta r) = C(\Delta r) = 10\% \cdot |\Delta r| +0.1\% \cdot \mathbf 1_{\Delta r\neq 0} $. The set of all central bank rates is $\mathcal R = (r^{c,i})_{i=1}^3 = [2.5\%,3\%,3.5\%]$, the transition function of the central bank rate $r^c$ is given by $P^c(r_{t+\Delta t}^c=r^{c,i'}|r_t^{c}=r^{c,i})=\lambda \Delta t \mathbf{1}_{i=i'} + (1-\lambda \Delta t)\mathbf{1}_{i\neq i'} $ with $\lambda = 0.2$, and initially we fix $r^c_i = 3\%$.

For the neural network architecture, we deploy a fully connected network with one hidden layer and a width of $L=256$. We set the outer iteration number to $N=100$, the inner iteration number to $M=400$, and the batch size to $B = 240$. The Adam algorithm is used for training, with learning rates set to $0.001$. In our experiment, the training process takes 12 minutes to complete on the CPU Intel Core i7-9750H.

\subsection{Experiment Result}
Under our parameter settings, we observed a decrease in the Bellman loss on batch, as shown in Figure \ref{fig:error-gap}. The Bellman loss for both the major player's Q-function and the minor player's functions alternates in descent but generally exhibits a decreasing trend. After about half of the training, the average loss remains below \(1 \times 10^{-7}\).

We also observe from Figures \ref{fig:dynamic_evolution} and \ref{fig:mean_field_evolution} that all banks consistently set lower deposit interest rates, regardless of fluctuations in the central bank policy rate. Under these parameter settings, banks prefer to minimize interest payments to depositors rather than compete for a larger market share. Especially, it can be observed from the approximated Nash equilibrium that all banks maintain their market share and benefit from the reduced interest expense. Moreover, Figure \ref{fig:mean_field_evolution}
shows a stable proportional movement, which can be attributed to the absence of interest rate differences.

Our major-minor mean-field game in bank interest rates involves a continuous state space and action space, bringing challenges to the projected mean-field approach proposed in \cite{cui2024learning}. This difficulty arises because an appropriate \textit{$\delta-$partition} would require a significant number of points to divide the probability space. In contrast, our neural network-based DQN algorithm effectively represents the Q functions by approximating the probability measure space by a 96-dimensional vector. Additionally, we propose a novel method for averaging the outputs of the neural networks within the Fictitious Play operation. Our algorithm has demonstrated stability in convergence and robustness in repeated tests, Bellman losses between iteration times $10,000$ and $20,000$, its maximum standard error over 10 tests is $7.23\times 10^{-7}$ (major player) and $2.48\times 10^{-7}$ (minor player) after $10000$ iterations. 

\section{Conclusion}
\label{sec:conclusion}

In this paper, we model the bank interest rate decision problem as an impulsive major-minor mean-field game against other players in the market. Banks adjust their interest rates to compete for market share, taking into account client attraction and deposit income. The banks are not symmetric, with finite major banks and multiple minor banks operating simultaneously. We developed a DQN algorithm to solve this major-minor mean-field game, introducing neural network approximation and averaging of networks. Our results demonstrate an approximated decision process of players under the Nash Equilibrium, providing insights into banks' decisions in adjusting interest rates to cope with competition.\footnote{\textbf{Disclaimer.} The authors' views are their own.
This paper was prepared for information purposes and is not a product of any Research Department. 
The authors' employers past and present make no representation and warranty whatsoever and disclaim all liability, for the completeness, accuracy, or reliability of the information contained herein. 
This document is not intended as investment research or investment advice, or a recommendation, offer or solicitation for the purchase or sale of any security, financial instrument, financial product or service, or to be used in any way to evaluate the merits of participating in any transaction, and shall not constitute a solicitation under any jurisdiction or to any person, if such solicitation under such jurisdiction or to such person would be unlawful.}

\bibliographystyle{ACM-Reference-Format}
\bibliography{ref}
\end{document}